\documentclass[12pt,twoside]{article} 
\usepackage{amsmath}
\usepackage{amsfonts}
\usepackage{amssymb}
\usepackage{graphicx}
\usepackage{tikz}

\tikzset{
	double color fill/.code 2 args={
		\pgfdeclareverticalshading[%
		tikz@axis@top,tikz@axis@middle,tikz@axis@bottom%
		]{diagonalfill}{100bp}{%
			color(0bp)=(tikz@axis@bottom);
			color(50bp)=(tikz@axis@bottom);
			color(50bp)=(tikz@axis@middle);
			color(50bp)=(tikz@axis@top);
			color(100bp)=(tikz@axis@top)
		}
		\tikzset{shade, left color=#1, right color=#2, shading=diagonalfill}
	}
}

\usetikzlibrary{patterns}
\tikzset{mynode/.style={fill,circle,inner sep=2pt,outer sep=0pt}}
\usetikzlibrary{decorations.markings}%
\tikzstyle{vertex} = [circle, draw, inner sep=0pt, minimum size=6pt]%
\usepackage{float}
\newcommand{\vertex}{\node[vertex]}%

\newcommand{\lb}{\linebreak}

\date{} 
\newcommand{\pf}{\noindent \textit{Proof}:\ } 
\newcommand{\qed}{\hfill{\ensuremath\blacksquare}}
\begin{document} 	
	\vspace*{-2.5cm}
	\centerline{\Large{\bf Clustering Coefficient of the Tensor }} 
	
	\centerline{\Large{\bf Product of Graphs}} 
	
	\centerline{\bf }
	
	\centerline{\bf {Remarl Joseph M. Damalerio and Rolito G. Eballe}} 
	
	\centerline{} 
	
	\centerline{Department of Mathematics} 
	
	\centerline{College of Arts and Sciences} 
	
	\centerline{Central Mindanao University} 
	
	\centerline{Musuan, Maramag, Bukidnon, 8714 Philippines}
	
	\newtheorem{Theorem}{\quad Theorem}[section] 
	
	\newtheorem{Definition}[Theorem]{\quad Definition} 
	
	\newtheorem{Corollary}[Theorem]{\quad Corollary} 
	
	\newtheorem{Lemma}[Theorem]{\quad Lemma} 
	
	\newtheorem{Example}[Theorem]{\quad Example}

	\begin{abstract} 
		Clustering coefficient is one of the most useful indices in complex networks. However, graph theoretic properties of this metric have not been discussed much in the literature, especially in graphs resulting from some binary operations. In this paper we present some expressions for the clustering coefficient of the tensor product of arbitrary graphs, regular graphs, and strongly regular graphs. A Vizing-type upperbound and a sharp lower bound for the clustering coefficient of the tensor product of graphs are also given.
		
	\end{abstract} 
	
	
	\vspace*{-.2cm}
	\section{Introduction} Let $G$ be a simple undirected graph with vertex set $V(G)$ and edge set $E(G)$. Let $N_G(v) = \{u\in V(G): uv\in E(G) \}$ be the open neighborhood of a vertex $v\in V(G)$, $\deg_Gv$ the degree of $v$, and $t_G(v) = |E(\langle N_G(v) \rangle)|$ the number of triangles in $G$ which are incident to $v$.  The \textit{local clustering coefficient} of vertex $v$ in $G$, denoted by $Cc_v(G)$, is a measure that assesses the local triangle density in a vertex's neighborhood. This number $Cc_v(G)$ can be defined as
	\vspace*{-.5cm}
	\begin{equation}
	Cc_v(G) = \begin{cases}
	\quad 0,  &\text{if }\deg_Gv \leq 1, \\ \frac{t_G(v)}{\binom{\deg_Gv}{2}}, & \text{if } \deg_Gv \geq 2.
	\end{cases}
	\label{lcc}
	\end{equation}
	
	\noindent This formula can be traced as a unifying version between its treatment in \cite{takahashi} and \cite{lietal}. On the other hand, the \textit{global clustering coefficient} $Cc(G)$ of a graph $G$ with order $n$ is a measure that indicates the overall clustering of $G$, obtained by averaging the local clustering coefficients of all the vertices in $G$. That is,  
	
	\begin{equation}
	Cc(G) = \frac{1}{n} \sum_{\substack{v \in V(G) \\ \deg_Gv\geq 2}}  Cc_v(G)  = \dfrac{1}{n}\sum_{\substack{v \in V(G) \\ \deg_Gv\geq 2}} \dfrac{2 t_G(v)}{\deg_Gv(\deg_G v -1)}.
	\label{gcc}
	\end{equation}
	
	This measure was introduced in the field of social network analysis by Duncan J. Watts and Steven Strogatz \cite{col_dy} in 1998 to determine whether a graph is a "small-world network". Since then, several studies from various standpoints have also emerged. However, as far as we know, there are no investigations or studies on the clustering coefficients of graphs resulting from some binary operations, although a related study on finding the number of distinct triangles in the tensor product $G\times H$ was done in \cite{trindex} while a triangle-counting algorithm for large networks appeared in \cite{suri2011counting}. \\
	\indent In this paper, we investigate the clustering coefficient of the tensor product of arbitrary graphs, regular graphs, and strongly regular graphs using some properties that the tensor product possesses and some inherent characteristics possessed by the factors or constituents. A Vizing-type upperbound and a sharp lower bound for the clustering coefficient of the tensor product of graphs are also aimed to be proved. Graphs considered in this paper are all finite and undirected simple graphs. For basic graph theory terminologies not specifically described nor defined in this paper, please refer to either \cite{chartrand1996graphs} or \cite{younger1972graph}.  
	
	\section{Tensor Product of Arbitrary Graphs}
	The \textit{tensor product} $G \times H$ of two graphs $G$ and $H$ is the graph with vertex set $V (G \times H) = V (G) \times V (H)$ and edge set $E(G \times H)$ satisfying the following adjacency condition: $(u,v)(u' ,v' ) \in E(G \times H)$ if and only if $uu' \in E(G)$ and $vv' \in E (H)$.
	\begin{Lemma} Let $G$ and $H$ be any graphs. If $u\in V(G)$ and $v\in V(H)$, then the number of triangles in $G\times H$ that are incident to the vertex $(u,v)\in V(G\times H)$ is given by the formula $ t_{G\times H}(u,v) = 2t_G(u)t_H(v).$
		\label{t(uv)}
	\end{Lemma}

	\pf Let $(u,v) \in V(G\times H)$ such that $(u',v')(u'',v'') \in E(\langle N_{G\times H}(u,v)\rangle )$. Then vertices $(u,v), (u',v'), (u'', v'') \in V(G\times H)$ are pairwise adjacent in $G\times H$. By the adjacency condition of the vertices in the tensor product, one can see that vertices $u,u',u''$ are pairwise adjacent in $G$ and so are the vertices $v, v', v''$ in $H$. Therefore, every triangle incident to $(u,v)$ in $G\times H$ emanates from a pair of distinct triangles, one of which is incident to $u$ in $G$ and another one incident to $v$ in $H$. \\
	\indent Now assume first that vertices $u, u', u''$ are pairwise adjacent in $G$ and $v,v',v''$ are pairwise adjacent in $H$. By the same adjacency condition in the tensor product we can see that vertices $(u',v'), (u'',v''), (u',v''),(u'',v') \in N_{G\times H}(u,v)$ such that  $(u',v')(u'',v''), (u',v'')(u'',v') \in E(\langle N_{G\times H}(u,v)\rangle)$. \lb Thus, each pair of triangles, one incident to $u$ in $G$ and the other incident to $v$ in $H$, produces two distinct triangles incident to $(u,v)$ in $G\times H$. That is, $\frac{1}{2}t_{G\times H}(u,v) = t_G(u)t_H(v)$. \qed \\
	
	In Figure 1 below, the highlighted vertices and edges  emphasized the two triangles incident to vertex $(u,v)$ in $G\times H$ as argued in Lemma \ref{t(uv)}.
	
	\begin{figure}[H]
		\label{tensex}
		\[\begin{tikzpicture}[mynode/.style={circle, inner sep=1pt,minimum size=6pt}, my node/.style={circle, inner sep=1pt,minimum size=6pt, double color fill={blue!80}{red!80}}]

		\node[my node, shading angle=45] () at (0,0) {};
		
		\node[] () at (-1.7,3) {$G \times H:$};
		\vertex (1) at (0,0) [label=left:$(u{,}v)$] {};
		\vertex (2) at (0,1) [label=left:$(u{,}v')$] {};
		\vertex (3) at (0,2) [label=left:$(u{,}v'')$] {};
		
		\vertex (4) at (2.3,0) [label=right:$(u'{,}v)$] {};
		\vertex (5) at (2.3,1) [label=right:$(u'{,}v')$] {};
		\vertex (6) at (2.3,2) [label=right:$(u'{,}v'')$] {};
		
		\vertex (7) at (4.6,0) [label=right:$(u''{,}v)$] {};
		\vertex (8) at (4.6,1) [label=right:$(u''{,}v')$] {};
		\vertex (9) at (4.6,2) [label=right:$(u''{,}v'')$] {};
		
		\draw[thick, dashed] (2.3,1) ellipse (4.9cm and 2.1cm);
		
		\draw[thick, dashed] (-5,1.95) ellipse (1.8cm and .9cm);
		\draw[thick, dashed] (-5,-.25) ellipse (1.8cm and .9cm);
		
		\node[mynode, fill=blue!80] () at (2.3,1) {}; 
		\node[mynode, fill=blue!80] () at (4.6,2) {};
		
		\node[mynode, fill=red!80] () at (2.3,2) {};
		\node[mynode, fill=red!80] () at (4.6,1) {};
		
		\node[] () at (-6.7,3) {$G :$};
		\vertex (10) at (-5,2.5) [label=right:$u$] {};
		\vertex (11) at (-4.3,1.7) [label=right:$u'$] {};
		\vertex (12) at (-5.7,1.7) [label=left:$u''$] {};
		
		\node[] () at (-6.7,.6) {$H :$};
		\vertex (13) at (-5,.3) [label=right:$v$] {};
		\vertex (14) at (-4.3,-.5) [label=right:$v'$] {};
		\vertex (15) at (-5.7,-.5) [label=left:$v''$] {};

		\path[thick]
		(2) edge[] (4) edge (6) edge (7) edge[] (9)
		(3) edge (4) edge (5) edge[bend right=50] (7) edge (8)
		(4) edge (8) edge (9)
		(5) edge (7) edge[blue,very thick] (9)
		(6) edge (7) edge[red,very thick] (8)
		
		(1) edge[blue,very thick] (5) edge[red!80,very thick] (6) edge[red,very thick] (8) edge[blue,very thick, bend left=50] (9)

		(10) edge (12) edge (11)
		(11) edge (12)
		(13) edge (14) edge (15)
		(15) edge (14)
		;
		\end{tikzpicture}\]
		\vspace*{-5mm}
		\caption{The tensor product $G \times H$ of two arbitrary graphs $G$ and $H$, in which vertices $u$, $u'$, $u''$ are pairwise adjacent in $G$ and vertices $v$, $v'$, $v''$ are pairwise adjacent in $H$.}
	\end{figure}
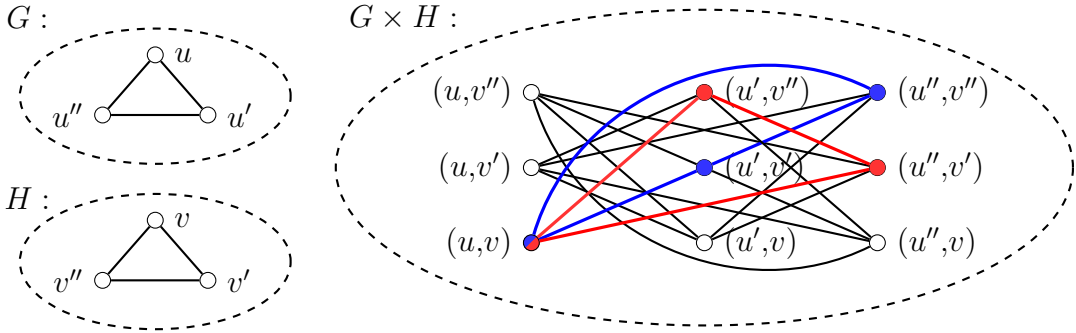 
	 
	\begin{Theorem} \label{lccGxH}
		Let $G$ and $H$ be any graphs. If $u\in V(G)$ and $v\in V(H)$ such that $\deg_Gu \geq 2$ and $\deg_Hv \geq 2$, then the local clustering coefficient of $(u,v)$ in $G \times H$ is given by the formula $$Cc_{(u,v)}(G\times H) = f(u,v) \cdot Cc_u(G) \cdot Cc_v(H),$$
		where $f(u,v) = (\deg_Gu-1)(\deg_Hv-1)/(\deg_Gu \cdot \deg_Hv - 1)$.
	\end{Theorem}
	\pf Using Equation (\ref{lcc}), Lemma \ref{t(uv)}, and the fact that $\deg_{G\times H}(u,v) = \deg_Gu \cdot \deg_Hv$, we have
	\begin{align*}
	Cc_{(u,v)}(G\times H) &= \frac{t_{G\times H}(u,v)}{\binom{\deg_{G\times H}(u,v)}{2}} = \frac{2t_G(u)t_H(v)}{\binom{\deg_Gu \cdot \deg_Hv}{2}} \\
	&= 2 \cdot Cc_u(G) \cdot Cc_v(H) \cdot \frac{\binom{\deg_Gu}{2}\binom{\deg_Hv}{2}}{\binom{\deg_Gu \cdot \deg_Hv}{2}} \\
	&= Cc_u(G) \cdot Cc_v(H) \cdot \frac{(\deg_Gu-1)(\deg_Hv-1)}{\deg_Gu \cdot \deg_Hv - 1}.
	\end{align*}
	
	\vspace*{-1cm}
	\qed 
	\\
	
	The next result, which is for the global clustering coefficient of $G\times H$, is a consequence of Theorem \ref{lccGxH}.
	
	\begin{Corollary}
		Let $G$ and $H$ be graphs of orders $n_1$ and $n_2$, respectively. Suppose $\delta(G) \geq 2$ and $\delta(H) \geq 2$. Then the global clustering coefficient of $G \times H$ is given by
		$$Cc(G\times H) = \frac{1}{n_1n_2} \sum_{u\in V(G)} \sum_{v\in V(H)} f(u,v) \cdot Cc_u(G) \cdot Cc_v(H),$$
		where $f(u,v) = (\deg_Gu-1)(\deg_Hv-1)/(\deg_Gu\deg_Hv - 1)$.
		\label{CcGxH}
	\end{Corollary}

	\pf Using Equation (\ref{gcc}) and Theorem \ref{lccGxH}, we obtain
	\begin{align*}
	Cc(G\times H) =& \frac{1}{n_1n_2} \sum_{(u,v)\in V(G\times H)} Cc_{(u,v)}(G\times H) \\
	=& \frac{1}{n_1n_2} \sum_{u\in V(G)} \sum_{v\in V(H)} Cc_u(G) \cdot Cc_v(H) \cdot \frac{(\deg_Gu-1)(\deg_Hu-1)}{(\deg_Gu \cdot \deg_Hv-1)}.
	\end{align*}
	The claimed equality follows.
	\qed \\
	
	The next result provides an upperbound for the global clustering coefficient of the tensor product of graphs, a Vizing-type relationship, albeit on the upperbound, with some restrictions on the factors.
	
	\begin{Corollary} \label{ub}
		For graphs $G$ and $H$ with $\delta(G) \geq 2$ and $\delta(H) \geq 2$, $$Cc(G\times H) \leq Cc(G) \cdot Cc(H);$$
		equality holds if at least one of $G$ and $H$ is a triangle-free graph. 
	\end{Corollary}
	\pf Let $u\in V(G)$ and $v\in V(H)$, where $G$ and $H$ are of orders $n_1$ and $n_2$, respectively. From Theorem \ref{lccGxH}, $Cc_{(u,v)}(G\times H) = f(u,v) \cdot Cc_u(G) \cdot Cc_v(H)$, where $f(u,v) = (\deg_Gu-1)(\deg_Hv-1)/(\deg_Gu \cdot \deg_Hv -1)$. Since $\deg_Gu, \deg_Hv \geq 2$, it follows that $(\deg_Gu-1)(\deg_Hv-1) < (\deg_Gu\cdot \deg_Hv) -1$. This means that
	$$0<f(u,v) = \frac{(\deg_Gu-1)(\deg_Hv-1)}{\deg_Gu\cdot \deg_Hv -1} < 1.$$
	Assume that each of the graphs $G$ and $H$ has a triangle. That is, there exists specific vertices $x\in V(G)$ and $y\in V(H)$ such that $t_G(x)\geq 1$ and $t_H(y)\geq 1$, implying that $Cc_{x}(G)> 0$ and $Cc_y(H) > 0$. Hence, for such vertices, it follows from Corollary \ref{lccGxH} that $$Cc_{(x,y)}(G\times H) = f(x,y)\cdot Cc_x(G)\cdot Cc_y(H) < Cc_x(G)\cdot Cc_y(H).$$
	As for the other vertices $(a,b)$ not incident to any triangle in $G\times H$, if there are any, their local clustering coefficients are clearly zero. Using Corollary \ref{CcGxH} and the inequality $f(u,v) <1$ above, we obtain
	\begin{align*}
	Cc(G\times H) =& \frac{1}{n_1n_2} \sum_{u\in V(G)} \sum_{v\in V(H)} f(u,v) \cdot Cc_u(G) \cdot Cc_v(H) \\
	<& \frac{1}{n_1n_2} \sum_{u\in V(G)} \sum_{v\in V(H)}  Cc_u(G) \cdot Cc_v(H) = Cc(G)\cdot Cc(H)
	\end{align*}
	Hence, $Cc(G\times H) < Cc(G)\cdot Cc(H)$, if both $G$ and $H$ have triangles.\\
	\indent On the other hand, assume that at least one of $G$ and $H$ is a triangle-free graph. Since the tensor product is commutative, we can assume that graph $G$ is triangle-free. Thus, for every $u\in V(G)$, $Cc_u(G) = 0$ and, hence, $Cc(G) = 0$. From this, we can see that $Cc(G)\cdot Cc(H) = 0$. From Corollary \ref{CcGxH}, we also have
	\begin{align*}
	Cc(G\times H) =& \frac{1}{n_1n_2} \sum_{u\in V(G)} \sum_{v\in V(H)} f(u,v) \cdot Cc_u(G) \cdot Cc_v(H) = 0
	\end{align*}
	Therefore, the relationship $Cc(G\times H)=Cc(G) \cdot Cc(H) $ is assured if at least one of $G$ and $H$ is a triangle-free graph. \qed \\
	
	The last paragraph of Corollary \ref{ub} above can actually be shortened by invoking Corollary 3.6 in \cite{trindex}, which says that if either $G$ or $H$ is triangle-free, then so is $G\times H$.
	
	In the 2016 paper of Yusheng Li, et al.\cite{lietal},  one of the results there gives a nontrivial lower bound for the clustering coefficient of $G$ with $\delta(G) \geq 2$. If $\delta (N_G(u))$ is the minimum degree of the subgraph of $G$ induced by $N_G(u)$ and if $\sigma = \min_{v\in V(G)} \delta (N_G(v))$, then $Cc(G) \geq \sigma /(\hat{d}-1),$ where $\hat{d}$ is the average degree of $G$ \cite{lietal}. An analogous result for $G\times H$ is given below.
	
	\begin{Corollary}\label{ntlb}
		Let $G$ and $H$ be graphs with $\delta(G)\geq 2$ and $\delta(H) \geq 2$. If $\sigma_G=\min_{u\in V(G)}\delta(N_G(u))$ and $\sigma_H=\min_{v\in V(H)}\delta(N_H(v))$, then
		$$Cc(G\times H) \geq \frac{\sigma_G \cdot \sigma_H}{\hat{d_G} \cdot \hat{d_H}-1},$$
		where $\hat{d_G}$ and $\hat{d_H}$ are the average degrees of $G$ and $H$, respectively. 
	\end{Corollary}  
	\pf Observe that 
	\begin{align*}
	t_G(u) = |E(\langle N_G(u) \rangle)|= \frac{1}{2} \sum_{w\in N_G(u)} \deg_{\langle N_G(u) \rangle}w \geq \frac{1}{2} \sum_{w\in N_G(u)} \sigma_G = \frac{\deg_Gu \cdot \sigma_G}{2}. 
	\end{align*}
	A similar argument applies to $t_H(v) \geq \frac{\deg_Hv~\cdot ~\sigma_H}{2}$. Additional simplifications yield $Cc_u(G) = \frac{t_G(u)}{\binom{\deg_Gu}{2}} \geq \frac{\sigma_G}{\deg_Gu-1}$ and $Cc_v(H) = \frac{t_H(v)}{\binom{\deg_Hv}{2}} \geq \frac{\sigma_H}{\deg_Hv-1}$. Applying these inequalities to Corollary \ref{CcGxH}, we obtain 
	
	\begin{align*}
	Cc(G\times H) =&\frac{1}{n_1n_2} \sum_{u\in V(G)} \sum_{v\in V(H)} f(u,v) \cdot Cc_u(G) \cdot Cc_v(H) \\
	\geq& \frac{1}{n_1n_2} \sum_{u\in V(G)} \sum_{v\in V(H)} \frac{(\deg_Gu -1)(\deg_Hv-1)}{\deg_Gu \cdot \deg_Hv-1} \cdot \frac{\sigma_G}{{\deg_Gu}-1} \cdot \frac{\sigma_H}{{\deg_Hv}-1} \\
	= & \frac{1}{n_1n_2} \sum_{u\in V(G)} \sum_{v \in V(G)} \frac{\sigma_G \cdot \sigma_H}{\deg_Gu \cdot \deg_Hv -1} \geq \frac{\sigma_G \cdot \sigma_H}{\hat{d_G}\cdot \hat{d_H} -1}. 
	\end{align*}

	The last inequality is due to the fact that $\deg_Gu \cdot \deg_Hv > 1$ for any $u \in V(G)$, $v\in V(H)$, and because the function $\frac{1}{x-1}$ is convex for $x>1$. \qed \\

\section{On Regular and Strongly Regular Graphs}
	A \textit{regular graph} is a graph that has uniform degree in its vertices. If $G$ is a regular graph with degree $d$ in all its vertices, then we call $G$ a $d$-regular graph. For general graphs, it is not viable to express the global clustering coefficient of the product cleanly in terms of the global clustering coefficients of its factors. But this is not the case for regular graphs as the next result shows.

	\begin{Theorem}\label{ccrg}
	Let $G$ and $H$ be graphs with orders $n_1$ and $n_2$, respectively. If $G$ and $H$ are regular graphs with respective degree regularity $d_G\geq 2$ and $d_H \geq 2$, then $$Cc(G\times H) = f \cdot Cc(G) \cdot Cc(H),$$ where $f = (d_G-1)(d_H-1)/(d_G\cdot d_H - 1)$.
	\end{Theorem}

	\pf Given that $G$ is $d_G$-regular and $H$ is $d_H$-regular, Corollary \ref{CcGxH} asserts that 
	\begin{align*}
	Cc(G\times H) =& \frac{1}{n_1n_2} \sum_{u\in V(G)} \sum_{v\in V(H)} \frac{(d_G-1)(d_H-1)}{(d_G\cdot d_H-1)} Cc_u(G) Cc_v(H) \\
	=& \frac{(d_G-1)(d_H-1)}{(d_G\cdot d_H-1)} \cdot \frac{1}{n_1} \sum_{u\in V(G)}   Cc_u(G) \cdot \frac{1}{n_2}\sum_{v\in V(H)} Cc_v(H) \\
	=& \frac{(d_G-1)(d_H-1)}{(d_G\cdot d_H-1)} \cdot Cc(G) \cdot Cc(H), 
	\end{align*}
	which completes the proof.
	\qed \\
	
	A special type of regular graphs are the \textit{strongly regular graphs}. Accordingly, a graph $G$ of order $n$ is said to be strongly regular with parameters $n$, $d$, $\mu_1$ ,$\mu_2$ , denoted $srg(n,d,\mu_1 ,\mu_2 )$, if it is $d$-regular, and any pair of vertices has $\mu_1$ common neighbors if they are
	adjacent, and $\mu_2$ common neighbors otherwise \cite{algrph}.
	
	\begin{Theorem}\cite{lietal}\label{srg}
		For any $srg(n,d,\mu_1,\mu_2)$ graph $G$ with $d \geq 2$, $$Cc(G) = \frac{\mu_1}{d-1}.$$
	\end{Theorem}
	\pf See Theorem 1 of \cite{lietal}. \qed \\
	
	Our last result says that Corollary \ref{ntlb} is sharp for the tensor product of strongly regular graphs. 
	
	\begin{Corollary}
		For two graphs $srg(n_1, d_G, \mu_1^G, \mu_2^G) ~ G$ and $srg(n_2, d_H, \mu_1^H, \mu_2^H) ~ H$ with $d_G,~ d_H \geq 2$,
		$$Cc(G\times H) = \frac{\mu_1^G \mu_1^H}{d_G\cdot d_H -1}.$$
	\end{Corollary}
	\pf Using Theorem \ref{ccrg} and Theorem \ref{srg}, we have 
	\begin{align*}
	Cc(G\times H) = \frac{(d_G-1)(d_H-1)}{d_G\cdot d_H-1}\cdot \frac{\mu_1^G}{d_G-1} \cdot \frac{\mu_1^H}{d_H-1} = \frac{\mu_1^G\mu_1^H}{d_G\cdot d_H-1},
	\end{align*} 
	where $\mu_1^G$ and $\mu_1^H$ correspond to $\sigma_G$ and $\sigma_H$ in Corollary \ref{ntlb}, respectively.\qed 
	
	\section{Final Remarks}  
	 Our motive in this work was to determine whether the parameter $Cc(G \times H)$ can be expressed meaningfuly in terms of $Cc(G)$ and $Cc(H)$, similar to our motives in \cite{acost} and \cite{trindex}. Our results in this paper showed that the global clustering coefficients of the factors are the key players, especially in the generated Vizing-type upperbound of $Cc(G\times H)$. For strongly regular graphs, a sharp lower bound was obtained. There are still a lot of work to be done in the clustering coefficient of graphs; we hope that this work could further stimulate research efforts into this area.
	\\

	{\bf Acknowledgement.} The authors would like to acknowledge the valuable comments and inputs made by the anonymous referee.

	{\bf Received: Month xx, 20xx}
	
\end{document}